\documentclass[11pt]{article}
\usepackage[english]{babel}
\usepackage{times}
\usepackage[T1]{fontenc}
\usepackage{amsmath, amsthm, amsfonts}
\usepackage{amssymb}
\usepackage{amscd}
\usepackage{xspace}
\usepackage{verbatim}
\usepackage{latexsym}
\usepackage{amsfonts}
\usepackage{graphicx}
\usepackage{float}
\usepackage{color}

\newcommand{\mysection}[1]{
\section{#1}\setcounter{equation}{0}}
\title{\bf The spherical $p$-harmonic eigenvalue problem \\in non-smooth domains}%%
\author{{\bf Konstantinos Gkikas\footnote{kgkikas@dim.uchile.cl}} \\[2mm]
 {\bf Laurent V\'eron\footnote{veronl@univ-tours.fr}}\\[2mm]
}%%ï¿½ï¿½Ä
\date{}
\begin{document}
 \maketitle
% \noindent{\small {\bf Abstract} We study the existence and uniqueness of  solutions of $\partial_tu-\Delta u+u^q=0$ ($q>1$) in $\Omega\times (0,\infty)$ where $\Omega\subset\mathbb R^N$ is a domain with a compact boundary, subject to the conditions $u=f\geq 0$ on $\partial\Omega\times (0,\infty)$ and the initial condition $\lim_{t\to 0}u(x,t)=\infty$. By means of Brezis' theory of maximal monotone operators in Hilbert spaces, we construct a minimal solution when $f=0$, whatever is the regularity of the boundary of the domain. When $\partial\Omega$ satisfies the parabolic Wiener criterion and $f$ is continuous, we construct a maximal solution and prove that it is the unique solution which blows-up at $t=0$.
% }

% \noindent
% {\it \footnotesize 1991 Mathematics Subject Classification}. {\scriptsize
% 35K60}.\\
% {\it \footnotesize Key words}. {\scriptsize Parabolic equations, singular solutions, semi-groups of contractions, maximal monotone operators, Wiener criterion.}
% \vspace{1mm}
% \hspace{.05in}

%% FONT commands

%% FONT commands
\newcommand{\txt}[1]{\;\text{ #1 }\;}%% Used in math only
\newcommand{\tbf}{\textbf}%% Bold face. Usage: \tbf{...}
\newcommand{\tit}{\textit}%% Italic
\newcommand{\tsc}{\textsc}%% Small caps
\newcommand{\trm}{\textrm}
\newcommand{\mbf}{\mathbf}%% Math bold
\newcommand{\mrm}{\mathrm}%% Math Roman
\newcommand{\bsym}{\boldsymbol}%% Bold math symbol
%%Macros for changing font size in math.
\newcommand{\scs}{\scriptstyle}%% as in subscript
\newcommand{\sss}{\scriptscriptstyle}%% as in sub-subscript
\newcommand{\txts}{\textstyle}
\newcommand{\dsps}{\displaystyle}
%%Macros for changing font size in text.
\newcommand{\fnz}{\footnotesize}
\newcommand{\scz}{\scriptsize}
%%\tiny<\scz<\fsz<\small<\large<\Large<\huge<\Huge
%%%%%%%%%%%%
%%%%%%%%%%%%
%% EQUATION commands
\newcommand{\be}{
\begin{equation}
}
\newcommand{\bel}[1]{
\begin{equation}
\label{#1}}
\newcommand{\ee}{
\end{equation}
}%% This macro does not work with amstex.
\newcommand{\eqnl}[2]{
\begin{equation}
\label{#1}{#2}
\end{equation}
}%%use not advisable; confusing
%%%%%%%%%%%%%%%
%% Unnumbered THEOREM env.
%% New env. to be used for unnumbered theorem, lemma etc. (but with specified name)
\newtheorem{subn}{\name}
\newcommand{\bsn}[1]{\def\name{#1}
\begin{subn}}
\newcommand{\esn}{
\end{subn}}
%%%%%%%%%%%%%%
%% NUMBERED THEOREM env.
%% Environments: theorem, lemma, corollary defintion and related commands,
%% designed to provide consecutive numbering of these forms.
\newtheorem{sub}{\name}[section]
\newcommand{\dn}[1]{\def\name{#1}}   %used in conjuction with sub or subn.
\newcommand{\bs}{
\begin{sub}}
\newcommand{\es}{
\end{sub}}
\newcommand{\bsl}[1]{
\begin{sub}\label{#1}}
%% the above must be preceeded by \dn (name definition),
%% however this is superceded by the list of commands bth etc.  below.
%%%%%%%%%%%%
%% NUMBERED THEOREM env. (cont.)
%% List of commands derived from 'sub' env. for theorem, lemma etc.
%% designed to provide consecutive numbering of these forms.
\newcommand{\bth}[1]{\def\name{Theorem}
\begin{sub}\label{t:#1}}
\newcommand{\blemma}[1]{\def\name{Lemma}
\begin{sub}\label{l:#1}}
\newcommand{\bcor}[1]{\def\name{Corollary}
\begin{sub}\label{c:#1}}
\newcommand{\bdef}[1]{\def\name{Definition}
\begin{sub}\label{d:#1}}
\newcommand{\bprop}[1]{\def\name{Proposition}
\begin{sub}\label{p:#1}}
%%%%%%%%%%%%%%%%%%%%%%%%%%%%%%%%%%
%% RERERENCE commands.
%% \newcommand{\R}[1]{$(\ref{#1})}
\newcommand{\R}{\eqref}
\newcommand{\rth}[1]{Theorem~\ref{t:#1}}
\newcommand{\rlemma}[1]{Lemma~\ref{l:#1}}
\newcommand{\rcor}[1]{Corollary~\ref{c:#1}}
\newcommand{\rdef}[1]{Definition~\ref{d:#1}}
\newcommand{\rprop}[1]{Proposition~\ref{p:#1}} 
%%%%%%%%%%%
%% ARRAY commands.
\newcommand{\BA}{
\begin{array}}
\newcommand{\EA}{
\end{array}}
\newcommand{\BAN}{\renewcommand{\arraystretch}{1.2}
\setlength{\arraycolsep}{2pt}
\begin{array}}
\newcommand{\BAV}[2]{\renewcommand{\arraystretch}{#1}
\setlength{\arraycolsep}{#2}
\begin{array}}
%Note: The first variable gives the amount of stretching: (#1) x default.
%For instance #1=1.2 means a 20% stretching. The second variable should be
%written for instance in the form  4pt ; here the default is 5pt
%\newcommand{\EAN}{\end{array}\setlength{\arraycolsep}{5pt}}
\newcommand{\BSA}{
\begin{subarray}}
\newcommand{\ESA}{\end{subarray}}
%Note: These are used in subscripts as well as superscripts. They work essentially
%% like 'array'.
\newcommand{\BAL}{\begin{aligned}}
\newcommand{\EAL}{\end{aligned}}
\newcommand{\BALG}{\begin{alignat}}
\newcommand{\EALG}{\end{alignat}}%% the abbrev. does not work with latex2e
\newcommand{\BALGN}{\begin{alignat*}}
\newcommand{\EALGN}{\end{alignat*}}%% the abbrev. does not work with latex2e
%% The 'aligned' environment must be placed inside an 'equation' env.
%% in the same way as the array.
%% One could use also the 'align' env. or the 'alignat' env.
%% However in this case each line is numbered, unless '\notag' is used.
%% The 'alignat'
%% has a slightly different format (the number of columns must be specified in advance)
%% but it has the advantage that the distance between columns is at our disposition.
%% (The default would be zero distance.) Using 'alignat*' we can have the advantages
%% of alignat plus the situation where separate lines are not numbered.
%% However in this case there is no numbering at all (unless we provide a tag).
%%%%%%%%%%
%% PROOF, REMARK etc.
\newcommand{\note}[1]{\textit{#1.}\hspace{2mm}}
\newcommand{\Proof}{\note{Proof}}
\newcommand{\qeda}{\hspace{10mm}\hfill $\square$}
\newcommand{\Remark}{\note{Remark}}
%%%%%%%% Style command.
\newcommand{\modin}{$\,$\\
[-4mm] \indent}
%% To be used after \mysection in order to start new line with \indent.
%%%%%%%%%%%%
%% MATHEMATICAL symbols
\newcommand{\forevery}{\quad \forall}
\newcommand{\set}[1]{\{#1\}}
\newcommand{\setdef}[2]{\{\,#1:\,#2\,\}}
\newcommand{\setm}[2]{\{\,#1\mid #2\,\}}
%% Arrows
\newcommand{\lra}{\longrightarrow}
\newcommand{\sgn}{\rm{sgn}}
\newcommand{\lla}{\longleftarrow}
\newcommand{\llra}{\longleftrightarrow}
\newcommand{\Lra}{\Longrightarrow}
\newcommand{\Lla}{\Longleftarrow}
\newcommand{\Llra}{\Longleftrightarrow}
\newcommand{\warrow}{\rightharpoonup}
%% Brackets, delimiters
\newcommand{
\paran}[1]{\left (#1 \right )}%% adjustable parantheses
\newcommand{\sqbr}[1]{\left [#1 \right ]}%% adjustable square brackets
\newcommand{\curlybr}[1]{\left \{#1 \right \}}%% adjustable curly brackets
\newcommand{\abs}[1]{\left |#1\right |}%% adjustable vertical delimiters
\newcommand{\norm}[1]{\left \|#1\right \|}%% adjustable norm
\newcommand{
\paranb}[1]{\big (#1 \big )}%% non-adjustable parantheses (big)
\newcommand{\lsqbrb}[1]{\big [#1 \big ]}%% non-adjustable square brackets (big)
\newcommand{\lcurlybrb}[1]{\big \{#1 \big \}}%% non-adjustable curly brackets (big)
\newcommand{\absb}[1]{\big |#1\big |}%% non-adjustable vertical delimiters (big)
\newcommand{\normb}[1]{\big \|#1\big \|}%% non-adjustable norm (big)
\newcommand{
\paranB}[1]{\Big (#1 \Big )}%% non-adjustable parantheses (Big)
\newcommand{\absB}[1]{\Big |#1\Big |}%% non-adjustable vertical delimiters (Big)
\newcommand{\normB}[1]{\Big \|#1\Big \|}%% non-adjustable norm (Big)

%%%%%%%%%%%%%%%%%
%% Adjustable parantheses etc. in a different DEFINITION format.
%\def\adp(#1){\left (#1 \right )}%% adjustable parantheses
%\def\adsb(#1){\left [#1\right ]}%% adjustable square brackets
%\def\adcb(#1){\left \{#1\right \}}%% adjustable curly brackets
%\def\abs|#1|{\left |#1\right |}%% adjustable vertical delimiters
%%%%%%%%%%%%%%%%
%% More mathematical symbols
\newcommand{\thkl}{\rule[-.5mm]{.3mm}{3mm}}
\newcommand{\thknorm}[1]{\thkl #1 \thkl\,}
\newcommand{\trinorm}[1]{|\!|\!| #1 |\!|\!|\,}
\newcommand{\bang}[1]{\langle #1 \rangle}%% angle bracket
\def\angb<#1>{\langle #1 \rangle}%% angle bracket
%% The two last lines yield the same result.
%% The second is used as follows: \angb<a,b>
\newcommand{\vstrut}[1]{\rule{0mm}{#1}}
\newcommand{\rec}[1]{\frac{1}{#1}}
%% OPERATOR names.
%% OPERATOR names.
\newcommand{\opname}[1]{\mbox{\rm #1}\,}
\newcommand{\supp}{\opname{supp}}
\newcommand{\dist}{\opname{dist}}
\newcommand{\myfrac}[2]{{\displaystyle \frac{#1}{#2} }}
\newcommand{\myint}[2]{{\displaystyle \int_{#1}^{#2}}}
\newcommand{\mysum}[2]{{\displaystyle \sum_{#1}^{#2}}}
\newcommand {\dint}{{\displaystyle \int\!\!\int}}%%%%%%%%%%
%%%%%%% SPACE commands
\newcommand{\q}{\quad}
\newcommand{\qq}{\qquad}
\newcommand{\hsp}[1]{\hspace{#1mm}}
\newcommand{\vsp}[1]{\vspace{#1mm}}
%%%%%%%%%%%
%% ABREVIATIONS
\newcommand{\ity}{\infty}
\newcommand{\prt}{
\partial}
\newcommand{\sms}{\setminus}
\newcommand{\ems}{\emptyset}
\newcommand{\ti}{\times}
\newcommand{\pr}{^\prime}
\newcommand{\ppr}{^{\prime\prime}}
\newcommand{\tl}{\tilde}
\newcommand{\sbs}{\subset}
\newcommand{\sbeq}{\subseteq}
\newcommand{\nind}{\noindent}
\newcommand{\ind}{\indent}
\newcommand{\ovl}{\overline}
\newcommand{\unl}{\underline}
\newcommand{\nin}{\not\in}
\newcommand{\pfrac}[2]{\genfrac{(}{)}{}{}{#1}{#2}}% frac with parantheses.
%%%%%%%%%%%
%%%%%%%%%%%%%

%%Macros for Greek letters.
\def\ga{\alpha}     \def\gb{\beta}       \def\gg{\gamma}
\def\gc{\chi}       \def\gd{\delta}      \def\ge{\epsilon}
\def\gth{\theta}                         \def\vge{\varepsilon}
\def\gf{\phi}       \def\vgf{\varphi}    \def\gh{\eta}
\def\gi{\iota}      \def\gk{\kappa}      \def\gl{\lambda}
\def\gm{\mu}        \def\gn{\nu}         \def\gp{\pi}
\def\vgp{\varpi}    \def\gr{\rho}        \def\vgr{\varrho}
\def\gs{\sigma}     \def\vgs{\varsigma}  \def\gt{\tau}
\def\gu{\upsilon}   \def\gv{\vartheta}   \def\gw{\omega}
\def\gx{\xi}        \def\gy{\psi}        \def\gz{\zeta}
\def\Gg{\Gamma}     \def\Gd{\Delta}      \def\Gf{\Phi}
\def\Gth{\Theta}
\def\Gl{\Lambda}    \def\Gs{\Sigma}      \def\Gp{\Pi}
\def\Gw{\Omega}     \def\Gx{\Xi}         \def\Gy{\Psi}

%%Macros for calligraphic letters.
\def\CS{{\mathcal S}}   \def\CM{{\mathcal M}}   \def\CN{{\mathcal N}}
\def\CR{{\mathcal R}}   \def\CO{{\mathcal O}}   \def\CP{{\mathcal P}}
\def\CA{{\mathcal A}}   \def\CB{{\mathcal B}}   \def\CC{{\mathcal C}}
\def\CD{{\mathcal D}}   \def\CE{{\mathcal E}}   \def\CF{{\mathcal F}}
\def\CG{{\mathcal G}}   \def\CH{{\mathcal H}}   \def\CI{{\mathcal I}}
\def\CJ{{\mathcal J}}   \def\CK{{\mathcal K}}   \def\CL{{\mathcal L}}
\def\CT{{\mathcal T}}   \def\CU{{\mathcal U}}   \def\CV{{\mathcal V}}
\def\CZ{{\mathcal Z}}   \def\CX{{\mathcal X}}   \def\CY{{\mathcal Y}}
\def\CW{{\mathcal W}} \def\CQ{{\mathcal Q}} 
%%%%%
%%Macros for 'blackboard' letters (See (27) for display.)
\def\BBA {\mathbb A}   \def\BBb {\mathbb B}    \def\BBC {\mathbb C}
\def\BBD {\mathbb D}   \def\BBE {\mathbb E}    \def\BBF {\mathbb F}
\def\BBG {\mathbb G}   \def\BBH {\mathbb H}    \def\BBI {\mathbb I}
\def\BBJ {\mathbb J}   \def\BBK {\mathbb K}    \def\BBL {\mathbb L}
\def\BBM {\mathbb M}   \def\BBN {\mathbb N}    \def\BBO {\mathbb O}
\def\BBP {\mathbb P}   \def\BBR {\mathbb R}    \def\BBS {\mathbb S}
\def\BBT {\mathbb T}   \def\BBU {\mathbb U}    \def\BBV {\mathbb V}
\def\BBW {\mathbb W}   \def\BBX {\mathbb X}    \def\BBY {\mathbb Y}
\def\BBZ {\mathbb Z}   \def\BBQ {\mathbb Q}

%%Macros for Ghotic (Fraktur) letters.
\def\GTA {\mathfrak A}   \def\GTB {\mathfrak B}    \def\GTC {\mathfrak C}
\def\GTD {\mathfrak D}   \def\GTE {\mathfrak E}    \def\GTF {\mathfrak F}
\def\GTG {\mathfrak G}   \def\GTH {\mathfrak H}    \def\GTI {\mathfrak I}
\def\GTJ {\mathfrak J}   \def\GTK {\mathfrak K}    \def\GTL {\mathfrak L}
\def\GTM {\mathfrak M}   \def\GTN {\mathfrak N}    \def\GTO {\mathfrak O}
\def\GTP {\mathfrak P}   \def\GTR {\mathfrak R}    \def\GTS {\mathfrak S}
\def\GTT {\mathfrak T}   \def\GTU {\mathfrak U}    \def\GTV {\mathfrak V}
\def\GTW {\mathfrak W}   \def\GTX {\mathfrak X}    \def\GTY {\mathfrak Y}
\def\GTZ {\mathfrak Z}   \def\GTQ {\mathfrak Q}

\font\Sym= msam10 % special symbols
\def\SYM#1{\hbox{\Sym #1}}
\newcommand{\bdw}{\prt\Gw\xspace}
\medskip
\abstract
{We prove the existence of p-harmonic functions under the form $u(r,\gs)=r^{-\gb}\gw(\gs)$ in any cone $C_S$ generated by a spherical domain $S$ and vanishing on $\prt C_S$. We prove the uniqueness of the exponent $\gb$ and of the normalized function $\gw$ under a Lipschitz condition on $S$.}\smallskip

\tableofcontents
\mysection{Introduction}

Let $p>1$, $S$ a domain of the unit sphere $S^{N-1}$ of $\BBR^N$ and $C_S:=\{(r,\gs):r>0,\gs\in S\}$  the positive cone generated by $S$. If one looks for $p$-harmonic functions in $C_S$ under the form $u(x)=u(r,\gs)=r^{-\gb}\gw(\gs)$ vanishing on 
$\prt C_S\setminus\{0\}$, then $\gw$ satisfies the {\it spherical $p$-harmonic eigenvalue problem} on $S$
\bel{I-1}\BA{lll}
-div'\left(\left(\gb^2\gw^2+\abs{\nabla'\gw}^2\right)^{\frac{p-2}{2}}\nabla'\gw\right)=(p-1)\gb(\gb-\gb_0)\left(\gb^2\gw^2+\abs{\nabla'\gw}^2\right)^{\frac{p-2}{2}}\gw\quad&\text{in }S\\
\phantom{-div'\left(\left(\gb^2\gw^2+\abs{\nabla'\gw}^2\right)^{\frac{p-2}{2}}\nabla'\right)}
\gw=0\quad&\text{in }\prt S
\EA\ee
with $\gb_0=\frac{N-p}{p-1}$ and were $div'$ and $\nabla'$ denote the divergence operator and the covariant gradient
on $S^{N-1}$ endowed with the metric induced by its isometric inbedding into $\BBR^N$. Separable solutions play a key role for describing the boundary behaviour  and the singularities of solutions of a large variety of quasilinear equations. When $N=2$ the equation is completely integrable and has been solved by Kroll in the regular case $\gb<0$ and Kichenassamy and V\'eron in the the singular case $\gb>0$. In higher dimension, Tolksdorff \cite{Tolk} proved the following:\smallskip

%%%%%%%%%%%%%%%%%%%%%%%%%%%%%%%%%%%%%%%%%%%%%%%%%%%%%
\nind {\bf Theorem A }{\it  If $S$ is a smooth spherical domain, there exist two couples $(\gb_S,\gw_S)$ and $(\gb'_S,\gw'_S)$ where 
$\gb_S>0$ and $\gb'_S<0$, $\gw_S$ and $\gw'_S$ are positive $C^2(\overline S)$-functions vanishing on $\prt S$
which solve $(\ref{I-1})$ with $(\gb,\gw)=(\gb_S,\gw_S)$ or $(\gb,\gw)=(\gb'_S,\gw'_S)$. Furthermore $\gb_S$ and $\gb'_S$ are unique, and $\gw_S$ and $\gw'_S$ are unique up to an homothety}.\smallskip

%%%%%%%%%%%%%%%%%%%%%%%%%%%%%%%%%%%%%%%%%%%%%%%%%%%%%
\nind A more general and transparent proof has been obtained by Porretta and V\'eron \cite{pv}, but always in the case of a smooth spherical domain. The aim of this article is to extend Theorem A to a general spherical domain. If we consider an increasing sequence of smooth domains $\{S_k\}$ such that $S_k\subset\overline S_k\subset S_{k+1}$ and  
$\cup_k S_k=S$ we prove the following:\smallskip
%%%%%%%%THEOREM %B%%%%%%%%%%%%%%%%%%%%%%%%%%%%%%%%%%%%%%%%%%%%%%%%%%%%%%%%%%%%%%%%%%%%%%%%%%%%%%%%%%%%%%%%%%%%

\nind {\bf Theorem B }{\it Assume that $S^c$ is not polar. Then the sequence of the $\gb_{S_k}>0$ from Theorem A is decreasing and converges to $\gb_S>0$. There exists $\gw_S\in W^{1,p}_0(S)\cap L^\infty(S)$ weak solution of $(\ref{I-1})$ with $\gb=\gb_S$. Furthermore $\gb_S>0$ is the largest exponent $\gb$ such that $(\ref{I-1})$ admits a positive solution $\gw_S\in W^{1,p}_0(S)$.}\smallskip 

%%%%%%%%THEOREM %C%%%%%%%%%%%%%%%%%%%%%%%%%%%%%%%%%%%%%%%%%%%%%%%%%%%%%%%%%%%%%%%%%%%%%%%%%%%%%%%%%%%%%%%%%%%%
Under a mild assumption on $S$ it is possible to approximate it by a decreasing sequence of smooth domains $S'_k$ such that $S'_k\subset \overline S'_k\subset S'_{k-1}$ and $\cap_k S'_k=\overline S$\smallskip

\nind {\bf Theorem C }{\it Assume that $S=\overset{o}{\overline S}$. Then the sequence $\gb_{S'_k}>0$ is increasing and converges to $\hat\gb_S>0$ and there exists $\hat\gw_S\in W^{1,p}_0(S)\cap L^\infty(S)$ weak solution of 
$(\ref{I-1})$ with $\gb=\hat\gb_S$. Furthermore $\hat\gb_S$ is the smallest exponent $\gb$ such that $(\ref{I-1})$ admits a positive solution $\gw_S\in W^{1,p}_0(S)$.}\smallskip 
%%%%%%%%%%%%%%%%%%%%%%%%%%%%%%%%%%%%%%%%%%%%%%%%%%%%%

We prove the uniqueness of the exponent $\gb$, under a Lipschitz assumption on $S$.  \smallskip

\nind {\bf Theorem D }{\it Assume that $S$ is a Lipschitz domain, then $\gb_S=\hat\gb_S$ and if $\gw$ and $\gw'$ are two  positive solutions of $(\ref{I-1})$ in $W^{1,p}_0(S)$, there exists a constant $c>0$ such that $c^{-1}\gw'\leq \gw\leq c\gw'$}. \smallskip

The proof of Theorem C is based upon a sharp form of boundary Harnack inequality proved in \cite {lip}, 
\bel{I-2}\BA{lll}
\abs{\ln\frac{\gw(\gs_1)}{\gw'(\gs_1)}-\ln\frac{\gw(\gs_2)}{\gw'(\gs_2)}}\leq c_1\abs{\gs_1-\gs_2}^\ga\quad\forall\, \gs_1,\gs_2\in S,
\EA\ee
for some $c_1=c_1(N,p,S)>0$ and $\ga\in (0,1)$. Actually we have a stronger result, much more delicate to obtain.\smallskip

\nind {\bf Theorem E }{\it Let $S$ be a Lipschitz subdomain of $S^{N-1}$. Then two positive solutions of $(\ref{I-1})$ in $W^{1,p}_0(S)$ are proportional.}\smallskip

The proof is based upon a non trivial adaptation of a series of deep results of Lewis and  Nystr\"om \cite {lip} concerning the $p$-Martin boundary of domains. \medskip

\nind{\bf Acknowledgements} This article has been prepared with the support of the collaboration
programs ECOS C14E08.

\mysection{Existence}

\subsection{Estimates}
Through this article we assume that $S^c$ is not polar, or equivalently that it has positive $c^{S^{N-1}}_{1,p}$-capacity. 
\blemma{est} Assume $p>1$. Then any solution $\gw\in W^{1,p}_0(S)$ of $(\ref{I-1})$ satisfies
\bel{II-0}
\displaystyle
\norm\gw_{C^{\gg}(S)}\leq c_1\norm\gw_{L^p(S)},
\ee
if $p>N-1$ where $\gg=1-\frac{N-1}{p}$ if $p>N-1$ and
\bel{II-1}
\displaystyle
\norm\gw_{L^\infty(S)}\leq c_1\norm\gw_{L^p(S)},
\ee
if $1<p\leq N-1$, where $c_1>0$ depends on $p$, $N$, $\gb$.
\es
\nind\Proof Multiplying the equation by $\gw$ and using H\"older's inequality, we derive
\bel{II-2}\BA{lll}
\displaystyle
(i)\quad&\myint{S}{}\left(\gb^2\gw^2+\abs{\nabla'\gw}^2\right)^{\frac{p}{2}} dS
\leq \left(\gb(p\gb-(p-1)\gb_0)\right)^{\frac{p}{2}}\myint{S}{}\abs\gw^p  dS\quad\!\text{if }\,p\geq 2,\\[4mm]
(ii)\quad&\myint{S}{}\left(\gb^2\gw^2+\abs{\nabla'\gw}^2\right)^{\frac{p}{2}} dS
\leq\gb^{p-1}(p\gb-(p-1)\gb_0)\myint{S}{}\abs\gw^p  dS\quad\text{if }\,1<p< 2.
\EA\ee
Notice that these inequalities hold for all $p>1$. If $p>N-1$ $(\ref{II-0})$ follows by Morrey'inequality. Here after we assume $1<p\leq N-1$. Let $\ga\geq 1$ and $k>0$. Then $\gz=\min\{\abs\gw,k\}^{\ga-1}\gw$ is an admissible test function, hence\smallskip

\nind1- If $p\geq 2$, 
\bel{II-3}\BA{lll}
\myint{S}{}\left(\gb^2\gw^2+\abs{\nabla'\gw}^2\right)^{\frac{p-2}{2}}\langle\nabla'\gw.\nabla'\gz\rangle dS=
(p-1)\gb(\gb-\gb_0)\myint{S}{}\left(\gb^2\gw^2+\abs{\nabla'\gw}^2\right)^{\frac{p-2}{2}}\gw\gz dS\\[4mm]
\phantom{---}
\leq c_2\myint{S}{}\abs{\nabla'\gw}^{p-2}\gw^2\min\{\abs\gw,k\}^{\ga-1} dS+c_2\gb^p\myint{S}{}\abs\gw^p\min\{\abs\gw,k\}^{\ga-1} dS\\[4mm]
\phantom{---}
\leq c_2\left(\myint{S}{}\abs\gw^p\min\{\abs\gw,k\}^{\ga-1} dS\right)^{\frac{p-2}{p}}\left(\myint{S}{}\abs{\nabla'\gw}^{p}\min\{\abs\gw,k\}^{\ga-1} dS\right)^{\frac{2}{p}}\\[4mm]
\phantom{---}
+ c_2\gb^p\myint{S}{}\abs\gw^p\min\{\abs\gw,k\}^{\ga-1} dS,
\EA\ee
where $c_2=c_2(N,p,\gb)>0$. Since
$$\myint{S}{}\left(\gb^2\gw^2+\abs{\nabla'\gw}^2\right)^{\frac{p-2}{2}}\langle\nabla'\gw.\nabla'\gz\rangle dS
\geq c_3(p)\myint{S}{}\abs{\nabla'\gw}^p\min\{\abs\gw,k\}^{\ga-1} dS,
$$
it implies that there exists $c_4=c_4(N,p,\gb)$ such that 
\bel{II-4}\BA{lll}
\myint{S}{}\abs{\nabla'\gw}^p\min\{\abs\gw,k\}^{\ga-1} dS\leq c_4\myint{S}{}\abs\gw^p\min\{\abs\gw,k\}^{\ga-1} dS,
\EA\ee
which yields
\bel{II-5}\BA{lll}
\myint{S}{}\abs{\nabla'j(\gw)}^p dS\leq c_4\myint{S}{}\abs{ j(\gw)}^pdS,
\EA\ee
where $j(\gw)=\min\{\abs\gw,k\}^{\frac{\ga-1}{p}}\gw$.\smallskip

\nind 2- If $1<p<2$, then
\bel{II-6}\BA{lll}
\myint{S}{}\left(\gb^2\gw^2+\abs{\nabla'\gw}^2\right)^{\frac{p-2}{2}}\langle\nabla'\gw.\nabla'\gz\rangle dS=
\myint{S}{}\left(\gb^2\gw^2+\abs{\nabla'\gw}^2\right)^{\frac{p-2}{2}}\abs{\nabla'\gw}^2 \min\{\abs\gw,k\}^{\ga-1}dS\\[4mm]\phantom{\myint{S}{}\left(\gb^2\gw^2+\abs{\nabla'\gw}^2\right)^{\frac{p-2}{2}}\nabla'\gw.}
+(\ga-1)\myint{S\cap\{\abs\gw<k\}}{}\left(\gb^2\gw^2+\abs{\nabla'\gw}^2\right)^{\frac{p-2}{2}}\abs{\nabla'\gw}^2 \abs\gw^{\ga-1} dS.
\EA\ee
Since
$$\BA{lll}
\myint{S}{}\left(\gb^2\gw^2+\abs{\nabla'\gw}^2\right)^{\frac{p-2}{2}}\abs{\nabla'\gw}^2 \min\{\abs\gw,k\}^{\ga-1}dS
=\myint{S}{}\left(\gb^2\gw^2+\abs{\nabla'\gw}^2\right)^{\frac{p}{2}}\min\{\abs\gw,k\}^{\ga-1}dS\\[4mm]\phantom{-------}
-\gb^2\myint{S}{}\left(\gb^2\gw^2+\abs{\nabla'\gw}^2\right)^{\frac{p-2}{2}} \min\{\abs\gw,k\}^{\ga-1}\gw^2dS
\\[4mm]\phantom{-----}
\geq \myint{S}{}\abs{\nabla'\gw}^p\min\{\abs\gw,k\}^{\ga-1} dS-\gb^2\myint{S}{}\left(\gb^2\gw^2+\abs{\nabla'\gw}^2\right)^{\frac{p-2}{2}} \min\{\abs\gw,k\}^{\ga-1}\gw^2dS,
\EA$$
we derive
\bel{II-7}\BA{lll}
\myint{S}{}\abs{\nabla'\gw}^p\min\{\abs\gw,k\}^{\ga-1} dS\leq \gb^{p-1}(p\gb-(p-1)\gb_0)
\myint{S}{}\abs{\gw}^p\min\{\abs\gw,k\}^{\ga-1}dS,
\EA\ee
which leads to $(\ref{II-5})$. Letting $k\to\infty$ we infer by Fatou's lemma,
\bel{II-8}\BA{lll}
\myint{S}{}\abs{\nabla'\abs{\gw}^{\frac{\ga-1}{p}+1}}^p dS\leq c_4\myint{S}{}\abs{\gw}^{\ga-1+p}dS.
\EA\ee
If $p<N-1$ we derive from Sobolev inequality and putting $q=\ga-1+p$ and $s=\frac{N-1}{N-1-p}>1$
\bel{II-9}\BA{lll}
\left(\myint{S}{}\abs{\gw}^{sq} dS\right)^{\frac1s}\leq c_5\myint{S}{}\abs{\gw}^{q}dS,
\EA\ee
and $c_5>0$ depends on $N$, $p$ and $\gb$. Iterating this estimate by Moser's method we derive $(\ref{II-9})$. \\
If $p=N-1$ we have for $1\leq m<p-1$ and $m^*=\frac{m(N-1)}{N-1-m}$
$$c_6\left(\myint{S}{}\abs{\gw}^{(\frac{\ga-1}{p}+1)m^*} dS\right)^{\frac{pm}{m^*}}\leq \left(\myint{S}{}\abs{\nabla'\abs{\gw}^{\frac{\ga-1}{p}+1}}^m dS\right)^{\frac{p}{m}}
\leq \abs S^{\frac{p}{m}-1}c_4\myint{S}{}\abs{\gw}^{\ga-1+p}dS,
$$
and $c_6=c_6(N,p)$, hence
\bel{II-10}\BA{lll}
\left(\myint{S}{}\abs{\gw}^{tq} dS\right)^{\frac1t}\leq c_5\myint{S}{}\abs{\gw}^{q}dS,
\EA\ee
with $t=\frac{m(N-1)}{p(N-1-m)}=\frac{m}{N-1-m}$. The proof follows again by Moser's iterative scheme.\qeda
%%%%%%%%%%%%%%%%%%%%%%%%%%%%%%%%%%%%%%%%%%%%%%%%%%%%%%%%%%%%%%%%%%%%%%%%%%%%%%%%%%%%%%%%%%%%%%%%%%%%%%%%%%%%%%%%%%%%%%%%%%%%%%%%%%%%%%%%%%%%%%%%%%%%%%%%%%%%%%%%%%%%%%%%%%%%%%%%%

\bprop{monot} Let $S_1$ and $S_2$ be two subdomains of $S^{N-1}$ such that $S_1\subset\overline S_1\subset S_2$ and $S_2$ not polar. Let $\gb_j>0$, j=1,2, such that there exist positive solutions $\gw_j\in W^{1,p}_0(S_j)$ solutions of 
\bel{II-11}\BA{lll}
-div'\left(\left(\gb_j^2\gw_j^2+\abs{\nabla'\gw_j}^2\right)^{\frac{p-2}{2}}\nabla'\gw_j\right)=(p-1)\gb_j(\gb_j-\gb_0)\left(\gb_j^2\gw_j^2+\abs{\nabla'\gw_j}^2\right)^{\frac{p-2}{2}}\gw_j\quad&\text{in }S_j\\
\phantom{-div'\left(\left(\gb^2\gw_j^2+\abs{\nabla'\gw_j}^2\right)^{\frac{p-2}{2}}\nabla'\right)}
\gw_j=0\quad&\text{in }\prt S_j.
\EA\ee
Then $\gb_1\geq\gb_2$.
\es
\nind\Proof Set $u_j(r,\gs)=r^{-\gb_j}\gw_j(\gs)$, $C_{S_j}=(0,\infty)\ti S_J$ and assume $\gb_1<\gb_2$. By Harnack inequality $\gw_2\geq c>0$ on 
$S_1$, thus 
$$u_2(r,\gs)\geq cr^{-\gb_2}\qquad\text{a.e. in }C_{S_1}.
$$
For $\ge>0$ there exist $r_\ge>0$ such that 
$$\ge u_2(x)\geq u_1(x)\qquad\forall\, x\in C_{S_1}\cap \overline B_{r_\ge}.
$$
Let $\gd>0$, there exists $R_\gd>0$ such that 
 $$ u_1(x)\leq\gd\qquad\forall\, x\in C_{S_1}\cap  B^c_{R_\gd}.
$$
Hence $\gz=(u_1-\ge u_2-\gd)_+\in W^{1,p}_0(Q^{r_\ge,R_\gd}_{S_1})$, where 
$Q^{r_\ge,R_\gd}_{S_1}=\{x\in C_{S_1}:r_\ge<\abs x<R_\gd\}$. This implies
$$\BA {lll}0=\myint{Q^{r_\ge,R_\gd}_{S_1}}{}\big\langle\abs{\nabla u_1}^{p-2}\nabla u_1-\abs{\nabla (\ge u_1)}^{p-2}\nabla (\ge u_1).\nabla\gz\big\rangle dx\\[4mm]
\phantom{0}
=\myint{Q^{r_\ge,R_\gd}_{S_1}\cap\{u_1-\ge u_2\geq\gd\}}{}
\big\langle\abs{\nabla u_1}^{p-2}\nabla u_1-\abs{\nabla (\ge u_1)}^{p-2}\nabla (\ge u_1).\nabla (u_1-u_2)\big\rangle dx.
\EA$$
Therefore $\nabla (u_1-\ge u_2-\gd)_+= 0$ a.e. in $Q^{r_\ge,R\gd}_{S_1}$, which leads to $u_1-\ge u_2\leq \gd$ in the same set. Letting $\gd\to 0$ yields $R_\gd\to\infty$, thus we obtain $u_1\leq \ge u_2$ in $C_{S_1}\cap \overline B^c_{r_\ge}$ hence $u_1\leq 0$ in
 $C_{S_1}$, contradiction.\qeda  %%%%%%%%%%%%%%%%%%%%%%%%%%%%%%%%%%%%%%%%%%%%%%%%%%%%%%%%%%%%%%%%%%%%%%%%%%%%%%%%%%%%%%%%%%%%%%%%%%%%%%%%%%%%%%%%%%%%%%%%%%%%%%%%%%%%%%%%%%%%%%%%%%%%%%%%%%%%%%%
 \subsection{Approximations from inside}
 \nind{\it Proof of Theorem B}. Let $\{S_k\}$ be an increasing sequence of smooth domains such that 
 $S_k\subset\overline S_k\subset S_{k+1}$. We denote by $\{(\gb_{S_k},\gw_k)\}$ the corresponding sequence of solutions of $(\ref{I-1})$ with $\gb=\gb_{S_k}$ and $\gw=\gw_k$. The sequence $\{\gb_{S_k}\}$ is uniquely determined by \cite{Tolk}, it admits a limit $\gb:=\gb_S$, and the $\gw_k$ are the unique positive solutions such that 
 $$\myint{S_k}{}\abs{\gw_k}dS=1.
 $$
 \nind If $p\geq 2$, we have
 $$\BA {lll}
 \myint{S_k}{}\abs{\nabla'\gw_k}^pdS\leq \myint{S_k}{}\left(\gb_{S_k}^2\gw_k^2+\abs{\nabla'\gw_k}^2\right)^{\frac{p-2}{2}}\abs{\nabla'\gw_k}^2dS
 \\[4mm]\phantom{ \myint{S_k}{}\abs{\nabla'\gw_k}^pdS}
=(p-1)\gb_{S_k}(\gb_{S_k}-\gb_0) \myint{S_k}{}\left(\gb_{S_k}^2\gw_k^2+\abs{\nabla'\gw_k}^2\right)^{\frac{p-2}{2}}\gw_k^2dS
 \\[4mm]\phantom{ \myint{S_k}{}\abs{\nabla'\gw_k}^pdS}
 \leq 2^{\frac{(p-4)_+}{2}}(p-1)\gb_{S_k}(\gb_{S_k}-\gb_0)\myint{S_k}{}\left(\gb_{S_k}^{p-2}\gw_k^p+\abs{\nabla'\gw_k}^{p-2}\gw_k^2\right)dS
  \\[4mm]\phantom{ \myint{S_k}{}\abs{\nabla'\gw_k}^pdS}
  \leq c_7(N,p,\gb_{S_k})\myint{S_k}{}\gw_k^p dS+\myfrac{1}{2}\myint{S_k}{}\abs{\nabla'\gw_k}^pdS.
 \EA$$
 Since $\gb_{S_k}\leq \gb_1$, we derive 
 \bel{III-1}\BA{lll}
 \myint{S_k}{}\abs{\nabla'\gw_k}^pdS\leq c_8,
\EA\ee
from the normalization assumption with $c_8=2c_7(N,p,\gb_1)$. \smallskip

 \nind If $1<p< 2$, we have
  $$\BA {lll}
   \myint{S_k}{}\abs{\nabla'\gw_k}^pdS\leq \myint{S_k}{}\left(\gb_{S_k}^2\gw_k^2+\abs{\nabla'\gw_k}^2\right)^{\frac{p}{2}}dS
   \\[4mm]\phantom{ \myint{S_k}{}\abs{\nabla'\gw_k}^pdS}
   \leq \gb_{S_k}(p\gb_{S_k}+(p-1)\gb_0)\myint{S_k}{}\left(\gb_{S_k}^2\gw_k^2+\abs{\nabla'\gw_k}^2\right)^{\frac{p-2}{2}}\gw_k^2dS
      \\[4mm]\phantom{ \myint{S_k}{}\abs{\nabla'\gw_k}^pdS}
      \leq  \gb^{p-1}_k(p\gb_{S_k}+(p-1)\gb_0)\myint{S_k}{}\gw_k^pdS,
  \EA$$
  and we obtain $(\ref{III-1})$ with $c_8=\gb^{p-1}_1(p\gb_1+(p-1)\gb_0)$.\smallskip
  
%%%%  
  Next we extend $\gw_k$ by $0$ in $S_k^c$. Then there exists $\gw\in W^{1,p}_0(S)$ such that $\gw_k\rightharpoonup \gw$ weakly  in $W^{1,p}_0(S)$, up to  subsequence that we still denote $\{\gw_k\}$, and $\gw_k\to \gw$ in $L^p(S)$. \smallskip
  %%%%%%%%%%%%%%%%%%%%%%%%%%%%%%%%%%%%%%%%%%%%%%%%%%%%%%%%%%%%%%%%%%%%%%%%%%%%%%%%%%%%%%%%%%%%%%%%%%%%%%%%%%%%%%%%%%%%%%%%%%%%%%%%%%%%%%%%%%%%%%%%%%%%%%%%%%%%%%%%%%%%%%%%%%%%%%%
  
  \nind{\it Step 1: We claim that $\nabla'\gw_k$ converges to $\nabla'\gw$ locally in $L^p(S)$}. \\
  Let $a\in S$ and $r>0$ such that  $B_{4r}(a)\subset S$. Then for $k\geq k_0$, $\overline B_{2r}(a)\subset S_k$. Let 
  $\gz\in C_0^\infty(B_{2r}(a))$ such that $0\leq\gz\leq 1$, $\gz=1$ in $B_{r}(a)$. For test function we choose 
  $\eta_k=\gz(\gw-\gw_k)$, then
  $$\myint{S_k}{}\left(\gb_{S_k}^2\gw_k^2+\abs{\nabla'\gw_k}^2\right)^{\frac{p-2}{2}}\langle\nabla'\gw_k.\nabla'\eta_k\rangle dS=
  (p-1)\gb_{S_k}(\gb_{S_k}-\gb_0)\myint{S_k}{}\left(\gb_{S_k}^2\gw_k^2+\abs{\nabla'\gw_k}^2\right)^{\frac{p-2}{2}}\gw_k\eta_k dS.
  $$
By the above inequality, we have
  $$\BA {lll}
  \myint{B_{2r}(a)}{}\bigg\langle\left(\gb^2\gw^2+\abs{\nabla'\gw}^2\right)^{\frac{p-2}{2}}\nabla'\gw-
\left(\gb_{S_k}^2\gw_k^2+\abs{\nabla'\gw_k}^2\right)^{\frac{p-2}{2}}\nabla'\gw_k.\nabla'\eta_k\bigg\rangle dS\\[4mm]
\phantom{------------}
= \myint{B_{2r}(a)}{}\left(\gb^2\gw^2+\abs{\nabla'\gw}^2\right)^{\frac{p-2}{2}}\langle\nabla'\gw.\nabla'\eta_k\rangle dS\\[4mm]
\phantom{------------}
-(p-1)\gb_{S_k}(\gb_{S_k}-\gb_0)\myint{S_k}{}\left(\gb_{S_k}^2\gw_k^2+\abs{\nabla'\gw_k}^2\right)^{\frac{p-2}{2}}\gw_k\eta_k dS.
  \EA$$
Using the weak convergence of the gradient, we have
$$\displaystyle
\lim_{k\to\infty}\myint{B_{2r}(a)}{}\left(\gb^2\gw^2+\abs{\nabla'\gw}^2\right)^{\frac{p-2}{2}}\langle\nabla'\gw.\nabla'\eta_k\rangle dS=0.
$$
Since $\gw_k$ is uniformly bounded in $W^{1,p}_0(S)$ and $\gw_k\to \gw$ in $L^p(S)$, we have
$$\displaystyle
\lim_{k\to\infty}\myint{B_{2r}(a)}{}\left(\gb_{S_k}^2\gw_k^2+\abs{\nabla'\gw_k}^2\right)^{\frac{p-2}{2}}\gw_k\eta_k dS=0,
$$
and
$$\displaystyle
\lim_{k\to\infty}\myint{B_{2r}(a)}{}(\gw-\gw_k)\bigg\langle\left(\gb^2\gw^2+\abs{\nabla'\gw}^2\right)^{\frac{p-2}{2}}\nabla'\gw-
\left(\gb_{S_k}^2\gw_k^2+\abs{\nabla'\gw_k}^2\right)^{\frac{p-2}{2}}\nabla'\gw_k.\nabla'\gz\bigg\rangle dS=0.
$$
Combining the above relations we infer
\bel{III-2}\BA {lll}\displaystyle
\lim_{k\to\infty}\myint{B_{2r}(a)}{}\gz\bigg\langle\left(\gb^2\gw^2+\abs{\nabla'\gw}^2\right)^{\frac{p-2}{2}}\nabla'\gw-
\left(\gb_{S_k}^2\gw_k^2+\abs{\nabla'\gw_k}^2\right)^{\frac{p-2}{2}}\nabla'\gw_k.\nabla'(\gw-\gw_k)\bigg\rangle dS=0.
\EA
\ee
Next we write
\bel{III-3}\BA {lll}
\myint{B_{2r}(a)}{}\gz\bigg\langle\left(\gb^2\gw^2+\abs{\nabla'\gw}^2\right)^{\frac{p-2}{2}}\nabla'\gw-
\left(\gb_{S_k}^2\gw_k^2+\abs{\nabla'\gw_k}^2\right)^{\frac{p-2}{2}}\nabla'\gw_k.\nabla'(\gw-\gw_k)\bigg\rangle dS\\[4mm]
\phantom{------}
=\myfrac{1}{2}
\myint{B_{2r}(a)}{}\gz\left(\left(\gb^2\gw^2+\abs{\nabla'\gw}^2\right)^{\frac{p-2}{2}}+
\left(\gb_{S_k}^2\gw_k^2+\abs{\nabla'\gw_k}^2\right)^{\frac{p-2}{2}}\right)\abs{\nabla'(\gw-\gw_k)}^2 dS\\[4mm]
\phantom{------}
+\myfrac{1}{2}
\myint{B_{2r}(a)}{}\gz\left(\!\left(\gb^2\gw^2+\abs{\nabla'\gw}^2\right)^{\frac{p-2}{2}}-
\left(\gb_{S_k}^2\gw_k^2+\abs{\nabla'\gw_k}^2\right)^{\frac{p-2}{2}}\right)\\[2mm]
\phantom{----------------------}
\ti\left(\abs{\nabla'\gw}^2+\gb^2\gw^2-\gb_{S_k}^2\gw_k^2-\abs{\nabla'\gw_k}^2\right) dS
\\[4mm]\phantom{------}
-\myfrac{1}{2}
\myint{B_{2r}(a)}{}\gz\left(\!\left(\gb^2\gw^2+\abs{\nabla'\gw}^2\right)^{\frac{p-2}{2}}-
\left(\gb_{S_k}^2\gw_k^2+\abs{\nabla'\gw_k}^2\right)^{\frac{p-2}{2}}\right)
\left(\gb^2\gw^2-\gb_{S_k}^2\gw_k^2\right) dS.
\EA\ee

\nind If $p\geq 2$, we have from $(\ref{II-3})$, 
\bel{III-4}\BA {lll}
\myint{B_{2r}(a)}{}\gz\bigg\langle\left(\gb^2\gw^2+\abs{\nabla'\gw}^2\right)^{\frac{p-2}{2}}\nabla'\gw-
\left(\gb_{S_k}^2\gw_k^2+\abs{\nabla'\gw_k}^2\right)^{\frac{p-2}{2}}\nabla'\gw_k.\nabla'(\gw-\gw_k)\bigg\rangle dS\\[4mm]
\phantom{------}
\geq \myfrac{1}{2}
\myint{B_{2r}(a)}{}\gz\left(\abs{\nabla'\gw}^{p-2}+
\abs{\nabla'\gw_k}^{p-2}\right)\abs{\nabla'(\gw-\gw_k)}^2 dS\\[4mm]\phantom{------}
-\myfrac{1}{2}
\myint{B_{2r}(a)}{}\gz\left(\!\left(\gb^2\gw^2+\abs{\nabla'\gw}^2\right)^{\frac{p-2}{2}}-
\left(\gb_{S_k}^2\gw_k^2+\abs{\nabla'\gw_k}^2\right)^{\frac{p-2}{2}}\right)
\left(\gb^2\gw^2-\gb_{S_k}^2\gw_k^2\right) dS
\\[4mm]
\phantom{------}
\geq \min\{2^{-1}, 2^{2-p}\}\myint{B_{2r}(a)}{}\gz\abs{\nabla'(\gw-\gw_k)}^pdS\\\phantom{------}
-\myfrac{1}{2}
\myint{B_{2r}(a)}{}\gz\left(\!\left(\gb^2\gw^2+\abs{\nabla'\gw}^2\right)^{\frac{p-2}{2}}-
\left(\gb_{S_k}^2\gw_k^2+\abs{\nabla'\gw_k}^2\right)^{\frac{p-2}{2}}\right)
\left(\gb^2\gw^2-\gb_{S_k}^2\gw_k^2\right) dS.
\EA\ee
Since $\gw_k\to\gw$ in $L^p(S)$, $\gb_{S_k}\to\gb$ and $\gw_k,\gw$ are uniformly bounded in $W^{1,p}_0(S)$ , we derive
$$\myint{B_{2r}(a)}{}\gz\left(\!\left(\gb^2\gw^2+\abs{\nabla'\gw}^2\right)^{\frac{p-2}{2}}-
\left(\gb_{S_k}^2\gw_k^2+\abs{\nabla'\gw_k}^2\right)^{\frac{p-2}{2}}\right)
\left(\gb^2\gw^2-\gb_{S_k}^2\gw_k^2\right) dS\to 0
$$
as $k\to\infty$. Jointly with $(\ref{III-2})$ we infer that 
\bel{III-5}\BA {lll}\displaystyle
\lim_{k\to\infty}\myint{B_r(a)}{}\abs{\nabla'(\gw-\gw_k)}^pdS=0.
\EA\ee
%%%%%%%%%%%%%%%%%%%%%%%%%%%%%%%%%%%%%%%%%%%%%%%%%%%%%%%%%%%%%%%%%%%%%%%%%%%%%%%%%%%%%%%%%%%%%%%%%%%%%%%%%%

\nind If $1<p< 2$, then
\bel{III-6}\BA {lll}\displaystyle
\myint{B_{2r}(a)}{}\gz\bigg\langle\left(\gb^2\gw^2+\abs{\nabla'\gw}^2\right)^{\frac{p-2}{2}}\nabla'\gw-
\left(\gb_{S_k}^2\gw_k^2+\abs{\nabla'\gw_k}^2\right)^{\frac{p-2}{2}}\nabla'\gw_k.\nabla'(\gw-\gw_k)\bigg\rangle dS\\[4mm]
\phantom{---}
=\myint{B_{2r}(a)}{}\gz\bigg\langle\left(\gb_{S_k}^2\gw_k^2+\abs{\nabla'\gw}^2\right)^{\frac{p-2}{2}}\nabla'\gw-
\left(\gb_{S_k}^2\gw_k^2+\abs{\nabla'\gw_k}^2\right)^{\frac{p-2}{2}}\nabla'\gw_k.\nabla'(\gw-\gw_k)\bigg\rangle dS\\[4mm]
\phantom{---}
+\myint{B_{2r}(a)}{}\gz\bigg\langle\left(\left(\gb^2\gw^2+\abs{\nabla'\gw}^2\right)^{\frac{p-2}{2}}-
\left(\gb_{S_k}^2\gw_k^2+\abs{\nabla'\gw}^2\right)^{\frac{p-2}{2}}\right)\nabla'\gw.\nabla'(\gw-\gw_k)\bigg\rangle dS.
\EA\ee
Up to extracting a subsequence, we have that $\gw_k\to\gw$ a.e. in $S$ and that there exists $\Gf\in L^1(S)$ such that 
 \bel{III-7}\BA {lll}\displaystyle
\abs{\gw_k}^p+\abs\gw^p\leq \Gf\qquad\text{a.e. in }S\quad\text{and }\;\forall\, k\geq 1.
\EA\ee
Since 
$$\left(\gb_{S_k}^2\gw_k^2+\abs{\nabla'\gw}^2\right)^{\frac{p-2}{2}}\abs{\nabla\gw}
\leq \left(\gb_{S_k}^2\gw_k^2+\abs{\nabla'\gw}^2\right)^{\frac{p-1}{2}}\leq \gb_{S_k}^{p-1}\gw^{p-1}_k+\abs{\nabla'\gw}^{p-1},
$$
and
$$\left(\gb^2\gw^2+\abs{\nabla'\gw}^2\right)^{\frac{p-2}{2}}\abs{\nabla\gw}\leq \gb^{p-1}\gw^{p-1}+\abs{\nabla'\gw}^{p-1},
$$
we derive that 
$$\abs{\left(\gb^2\gw^2+\abs{\nabla'\gw}^2\right)^{\frac{p-2}{2}}-
\left(\gb_{S_k}^2\gw_k^2+\abs{\nabla'\gw}^2\right)^{\frac{p-2}{2}}}\abs{\nabla'\gw}\leq 
2\left(\gb^{p-1}\Gf^{p-1}+\abs{\nabla'\gw}^{p-1}\right),
$$
which implies that 
$$\gz\left(\left(\gb^2\gw^2+\abs{\nabla'\gw}^2\right)^{\frac{p-2}{2}}-
\left(\gb_{S_k}^2\gw_k^2+\abs{\nabla'\gw}^2\right)^{\frac{p-2}{2}}\right)\nabla'\gw\to 0 \quad\text { in }L^{p'}(S)
$$
where $p'$ is the conjugate of $p$, and finally
\bel{III-8}\BA{ll}
\myint{B_{2r}(a)}{}\gz\langle\left(\left(\gb^2\gw^2+\abs{\nabla'\gw}^2\right)^{\frac{p-2}{2}}-
\left(\gb_{S_k}^2\gw_k^2+\abs{\nabla'\gw}^2\right)^{\frac{p-2}{2}}\right)\nabla'\gw.\nabla'(\gw-\gw_k)\rangle dS\to 0\;\;\text{as }k\to\infty.
\EA\ee
For the last term on the right-hand side of $(\ref{III-6})$, we have, for $\gg\in\BBR_+$ and ${\bf A},{\bf B}\in\BBR^N$,
$$\BA {lll}\left(\gg+\abs {\bf B}^2\right)^{\frac{p-2}{2}}\!\!{\bf B}-\left(\gg+\abs {\bf A}^2\right)^{\frac{p-2}{2}}\!\!{\bf A}
=\!\myint{0}{1}\myfrac{d}{dt}\left(\left(\gg+\abs {t{\bf B}+(1-t){\bf A}}^2\right)^{\frac{p-2}{2}}\!\!(t{\bf B}+(1-t){\bf A})\right) dt\\[4mm]
\phantom{\left(\gg+\abs {\bf B}^2\right)^{\frac{p-2}{2}}{\bf B}-\left(\gg+\abs {\bf A}^2\right)^{\frac{p-2}{2}}{\bf A}}
=\left(\myint{0}{1}\left(\gg+\abs {t{\bf B}+(1-t){\bf A}}^2\right)^{\frac{p-2}{2}}dt\right)({\bf B}-{\bf A})\\[4mm]
\phantom{\left(\gg+\abs {\bf B}\right)}
+(p-2)\myint{0}{1}\left(\gg+\abs {t{\bf B}+(1-t){\bf A}}^2\right)^{\frac{p-4}{2}}\langle t{\bf B}+(1-t){\bf A}\,.\,{\bf B}-{\bf A}\rangle(t{\bf B}+(1-t){\bf A})dt.
\EA$$
This implies
$$\BA {lll}
\!\!\langle\left(\gg+\abs {\bf B}^2\right)^{\frac{p-2}{2}}\!\!{\bf B}-\!\left(\gg+\abs {\bf A}^2\right)^{\frac{p-2}{2}}\!\!{\bf A}\,.\,{\bf B}-{\bf A}\rangle
=\left(\myint{0}{1}\!\!\left(\gg+\abs {t{\bf B}+(1-t){\bf A}}^2\right)^{\frac{p-2}{2}}\!\!\!\!dt\right)\abs{{\bf B}-{\bf A}}^2\\[4mm]
\phantom{\left(\abs {\bf }^2\right)^{\frac{p-2}{2}}{\bf B}-\left(\gg \right)^{\frac{p}{2}}{\bf A},}
+(p-2)\myint{0}{1}\!\!\left(\gg+\abs {t{\bf B}+(1-t){\bf A}}^2\right)^{\frac{p-4}{2}}\langle t{\bf B}+(1-t){\bf A}\,.\,{\bf B}-{\bf A}\rangle^2dt.
\EA$$
We observe that
$$\BA {lll}\myint{0}{1}\!\!\left(\gg+\abs {t{\bf B}+(1-t){\bf A}}^2\right)^{\frac{p-4}{2}}\langle t{\bf B}+(1-t){\bf A}\,.\,{\bf B}-{\bf A}\rangle^2dt
\\[4mm]
\phantom{-----------------}
\leq\abs{{\bf B}-{\bf A}}^2\myint{0}{1}\!\!\left(\gg+\abs {t{\bf B}+(1-t){\bf A}}^2\right)^{\frac{p-2}{2}}dt,
\EA$$
and since $1<p<2$, we finally obtain
\bel{III-9}\BA{ll}
\langle\left(\gg+\abs {\bf B}^2\right)^{\frac{p-2}{2}}\!\!{\bf B}-\!\left(\gg+\abs {\bf A}^2\right)^{\frac{p-2}{2}}\!\!{\bf A}\,.\,{\bf B}-{\bf A}\rangle\\[3mm]
\phantom{------------}
\geq (p-1)\left(\myint{0}{1}\left(\gg+\abs {t{\bf B}+(1-t){\bf A}}^2\right)^{\frac{p-2}{2}}dt\right)\abs{{\bf B}-{\bf A}}^2\\[4mm]
\phantom{------------}
\geq (p-1)\abs{{\bf B}-{\bf A}}^2\left(\gg+1+\abs {{\bf B}}^2+\abs{{\bf A}}^2\right)^{\frac{p-2}{2}}.
\EA\ee
We plug this estimate into $(\ref{III-6})$ with $\gg=\gb^2_k\gw^2_k$, ${\bf A}=\nabla'\gw$ and ${\bf B}=\nabla'\gw_k$, then
\bel{III-10}\BA{ll}
\myint{B_{2r}(a)}{}\gz\langle\left(\gb_{S_k}^2\gw_k^2+\abs{\nabla'\gw}^2\right)^{\frac{p-2}{2}}\nabla'\gw-
\left(\gb_{S_k}^2\gw_k^2+\abs{\nabla'\gw_k}^2\right)^{\frac{p-2}{2}}\nabla'\gw_k.\nabla'(\gw-\gw_k)\rangle dS\\[4mm]
\phantom{-------}
\geq \myint{B_{2r}(a)}{}\gz\abs{\nabla'(\gw-\gw_k)}^2\left(\gb^2_k\gw^2_k+1+\abs{\nabla'\gw_k}^2+\abs{\nabla'\gw}^2\right)^{\frac{p-2}{2}}dS.
\EA\ee
Set $\phi(.)=\gb^2_k\gw^2_k+1+\abs{\nabla'\gw_k}^2+\abs{\nabla'\gw}^2$, then
\bel{III-11}\BA{ll}
\myint{B_r(a)}{}\abs{\nabla'\gw-\nabla'\gw_k}^p dS=\myint{B_r(a)}{}\abs{\nabla'\gw-\nabla'\gw_k}^p \phi^{\frac{p(p-2)}{4}}\phi^{-\frac{p(p-2)}{4}}dS\\[4mm]
\phantom{\myint{B_r(a)}{}\abs{\nabla'\gw-\nabla'\gw_k}^p dS}
\leq \left(\myint{B_r(a)}{}\abs{\nabla'\gw-\nabla'\gw_k}^2 \phi^{\frac{p-2}{2}}dS\right)^{\frac{p}{2}}
\left(\myint{B_r(a)}{} \phi^{\frac{p}{2}}dS\right)^{\frac{2-p}{2}}.
\EA\ee
Jointly with $(\ref{III-2})$  and $(\ref{III-10})$ we conclude that  $(\ref{III-5})$. Step 1 follows by a standard covering argument.\medskip

 \nind{\it Step 2: We claim that $\gw_k$ converges to $\gw$ in $W^{1,p}_0(S)$}. \\
 %%%%%%%%
Up to a subsequence that we denote again by $\{k\}$, we can assume that $\gw_k\to\gw$ and $\nabla'\gw_k\to\nabla'\gw$ a.e. in $S$. Let $\gz\in C^\infty_0(S)$, then there exists $k_\ge\in\BBN$ such that the support $K$ of $\gz$ is a compact subset of 
$S_k$ for all $k\geq k_\ge$. If $1<p<2$,
$$\left(\gb_{S_k}^2\gw_k^2+\abs{\nabla'\gw_k}^2\right)^{\frac{p-2}{2}}\abs{\nabla'\gw_k}\leq \abs{\nabla'\gw_k}^{p-1},
$$
which bounded in $L^{p'}(K)$, then uniformly integrable in $K$ and by Vitali's convergence theorem
$$\left(\gb_{S_k}^2\gw_k^2+\abs{\nabla'\gw_k}^2\right)^{\frac{p-2}{2}}\nabla'\gw_k\to 
\left(\gb^2\gw^2+\abs{\nabla'\gw}^2\right)^{\frac{p-2}{2}}\nabla'\gw,$$
in $L^1_{loc}(S)$. Similarly
$$\left(\gb_{S_k}^2\gw_k^2+\abs{\nabla'\gw_k}^2\right)^{\frac{p-2}{2}}\gw_k\to 
\left(\gb^2\gw^2+\abs{\nabla'\gw}^2\right)^{\frac{p-2}{2}}\gw,$$
in $L^1_{loc}(S)$. If $p\geq 2$
$$\left(\gb_{S_k}^2\gw_k^2+\abs{\nabla'\gw_k}^2\right)^{\frac{p-2}{2}}\abs{\nabla'\gw_k}\leq c\left(\abs{\gw_k}^{p-1}+\abs{\nabla'\gw_k}^{p-1}\right),
$$
and we conclude again by Vitali's convergence theorem that the previous convergences hold. 
Since 
$$\myint{S_k}{}\left(\gb_{S_k}^2\gw_k^2+\abs{\nabla'\gw_k}^2\right)^{\frac{p-2}{2}}\langle\nabla'\gw_k.\nabla'\gz\rangle dS=
(p-1)\gb_{S_k}(\gb_{S_k}-\gb_0)\myint{S_k}{}\left(\gb_{S_k}^2\gw_k^2+\abs{\nabla'\gw_k}^2\right)^{\frac{p-2}{2}}\gw_k\gz dS
$$
we conclude that $\gw$ is a weak solution of $(\ref{I-1})$ with $\gb=\gb_S$.\qeda

 \subsection{Approximations from outside}
 \nind{\it Proof of Theorem C}. Since $\overline S^c$ has a non-empty interior, the existence of a sequence $\{\gw'_k\}$ corresponding to solutions of $(\ref{I-1})$ in $S'_k$ with $\gb=\gb_{S'_k}$ is the consequence of \cite{pv}. The fact that $\{\gb_{S'_k}\}$ is increasing follows from \rprop{monot}. We denote by $\hat\gb:=\hat \gb_S$ its limit, and it is smaller or equal to  $\gb_S$. Estimates $(\ref{II-3})$ are valid with $S'_k$, $\gw'_k$ and $\gb_{S'_k}$ instead of $S$, $\gw$ and $\gb$. If we extend $\gw'_k$ by $0$ in $S'^c_k$ these estimates are valid with $S^{N-1}$ instead of $S'_k$. Then up to a subsequence the exists $\gw\in W^{1,p}(S^{N-1})$ and a subsequence stil denoted by $\{k\}$ such that $\gw'_k\rightharpoonup\gw$ weakly in $W^{1,p}(S^{N-1})$, strongly in $L^{p}(S^{N-1})$ and a.e. in $S^{N-1}$. Furthermore, as in the proof of Theorem A, for any compact set $K\subset S$, $\nabla'\gw'_k\to\nabla'\gw'$ in $L^p(K)$. This is sufficient to assert that $\gw$ is a weak solution of 
 $$-div'\left(\left(\hat \gb^2\gw'^2+\abs{\nabla'\gw'}^2\right)^{\frac{p-2}{2}}\nabla'\gw'\right)=(p-1)\hat \gb(\hat \gb-\gb_0)\left(\hat \gb^2 \gw^2+\abs{\nabla'\gw'}^2\right)^{\frac{p-2}{2}}\gw'\quad\text{in }S.
 $$
 Moreover $\gw'\lfloor_{S'_k}$ belongs to $W_0^{1,p}(S'_k)$ for all $k$. Since $\gw'_k=0$ in $S_k^c$ and converges a.e. to $\gw$, this last function vanishes a.e. in $\cup_kS_k^c=(\cap_k S_k)^c=\overline S^c$. Therefore $\gw$ vanishes a.e. in $\overline S^c$ and since it is quasi continuous, it vanishes, $(1-p)$- quasi everywhere in $\overline S^c$. From Netrusov's theorem (see \cite[Th 10.1.1]{AdHed}-(iii)) there exists a sequence $\{\eta_n\}\subset C^\infty_0( S)$ which converges to $\gw$ in $W^{1,p}(S)$, thus $\gw\in W_0^{1,p}(S)$.\qeda
 
 %%%%%%%%%%%%%%%%%%%%%%%%%%%%%%%%%%%%%%%%%%%%%%%%%%%%%%%%%%%%%%%%%%%%%%%%%%%%%%%%%%%%%%%%%%%%%%%%%%%%%%%%%%%%%%
  \section{Uniqueness}
 %%%%%%%%%%%%%%%%%%%%%%%%%%%%%%%%%%%%%%%%%%%%%%%%%%%%%%%%%%%%%%%%%%%%%%%%%%%%%%%%%%%%%%%%%%%%%%%%%%%%%%%%%%%%%%%%%%%%%%%%%%%%%%
 \subsection{Uniqueness of exponent $\gb$}
 
  \nind{\it Proof of Theorem D}.  If $S$ is Lipschitz, $C_S$ is also Lipschitz. We fix $z\in S\approx S^{N-1}\cap \prt C_S$ and we apply \cite[Th 2]{lip} in $G_z=C_S\cap B_{\frac{1}{2}}(z)$ to two separable $p$-harmonic functions $u(r,\gs)=r^{-\gb}\gw(\gs)$ and $u'(r,\gs)=r^{-\gb'}\gw'(\gs)$. There exist $\gg\in (0,\frac{1}{2})$, $c_{10}>0$ and $\ga\in (0,1)$ such that 
  \bel{IV-1}\BA{ll}
\abs{\ln\myfrac{u(y_1)}{u'(y_1)}-\ln\myfrac{u(y_2)}{u'(y_2)}}\leq c_{10}\abs{y_1-y_2}^\ga\quad\forall\, y_1,y_2\in 
C_S\cap B_{\gg}(z).
\EA\ee
Assume $|y_1|=|y_2|=1$, then 
  \bel{IV-2}\BA{ll}
\abs{\ln\myfrac{\gw(y_1)}{\gw'(y_1)}-\ln\myfrac{\gw(y_2)}{\gw'(y_2)}}\leq c_{10}\abs{y_1-y_2}^\ga\quad\forall\, y_1,y_2\in 
S\cap B_{\gg}(z).
\EA\ee
We denote by $\ell (x,y)$ the geodesic distance on $S^{N-1}$ and by $\ell (x,K)$ the geodesic distance from a point $x\in S^{N-1}$ to a subset $K$.  Since the set $S_\gg=\{\gs\in S:\ell(\gs,\prt S)\leq\tfrac\gg 2\}$ can be covered by a finite number of balls with center on 
$\prt S$, we infer that 
  \bel{IV-3}\BA{ll}
\abs{\ln\myfrac{\gw(y_1)}{\gw'(y_1)}-\ln\myfrac{\gw(y_2)}{\gw'(y_2)}}\leq c_{11}\qquad\forall\, y_1,y_2\in 
S_\gg.
\EA\ee
In $S\setminus \overline S_{\tfrac\gg 2}$ we can use Harnack inequality to obtain
  \bel{IV-4}\BA{ll}-c_{12}\leq  \ln\myfrac{\gw(y_1)}{\gw(y_2)}\leq c_{12}\qquad\forall\, y_1,y_2\in S\setminus \overline S_{\tfrac\gg 2}\,\text{ s.t. }\ell(y_1,y_2)\leq \tfrac\gg 2.\EA\ee
Hence there exists a constant $c_{13}>0$ such that $(\ref{IV-4})$ holds for any $y_1,y_2\in S\setminus \overline S_{\tfrac\gg 2}$, with $c_{12}$ replaced by $c_{13}$. Furthermore $\gw'$ satisfies the same inequality in 
$S\setminus \overline S_{\tfrac\gg 2}$. Combining the two inequalities we obtain
  \bel{IV-5}\BA{ll}-2c_{13}\leq  \ln\myfrac{\gw(y_1)}{\gw(y_2)}-\ln\myfrac{\gw'(y_1)}{\gw'(y_2)}\leq 2c_{13}\qquad\forall\, y_1,y_2\in S\setminus \overline S_{\tfrac\gg 2}.\EA\ee
Combining this estimate with $(\ref{IV-2})$ we derive that it holds  for all $y_1,y_2\in S$. This implies 
  \bel{IV-6}\BA{ll}
  e^{-2c_{13}}\myfrac{\gw(y_2)}{\gw'(y_2)}\leq  \myfrac{\gw(y_1)}{\gw'(y_1)}\leq  e^{2c_{13}}\myfrac{\gw(y_2)}{\gw'(y_2)}\qquad\forall\, y_1,y_2\in S.\EA\ee
  
  Assume now that there exist two exponents $\gb>\gb'>0$ such that $r^{-\gb}\gw(.)$ and $r^{-\gb'}\gw'(.)$ are $p$-harmonic and positive in the cone $C_S$ and vanishes on $\prt C_S$. Put $\gth=\frac{\gb}{\gb'}$, $\eta=\gw'^\gth$ and 
 $$ \CT(\eta)=-div'\left(\left(\gb^2\eta^2+\abs{\nabla'\eta}^2\right)^{\frac{p-2}{2}}\nabla'\eta\right)-(p-1)\gb(\gb-\gb_0)\left(\gb^2\eta^2+\abs{\nabla'\eta}^2\right)^{\frac{p-2}{2}}\eta,
  $$
  then
  $$\CT(\eta)=-\gth^{p-2}\left(\gb'^2\gw'^2+\abs{\nabla'\gw'}^2\right)^{\frac{p-2}{2}}\left((\gb-\gb')\gw'^2+(p-1)\gth(\gth-1)\abs{\nabla'\gw'}^2\right)\leq 0.
  $$
  Up to multiplying $\gw'$ by $\gl$, we can assume that $\eta\leq\gw$ and that the graphs of $\eta$ and $\gw$ are tangent in $\overline S$. Since $\gw'\leq c\gw$, $\eta=o(\gw)$ near $\prt S$. Hence there exists $\gs_0\in S$ such that 
  $\gw(\gs_0)=\eta(\gs_0)$ and the coincidence set of $\eta$ and $\gw$ is a compact subset of $S$. We put $w=\gw-\eta$, since $\nabla\gw(\gs_0)=\nabla\eta(\gs_0)$ we proceed as in \cite[Th 4.1]{pv2} (see also \cite{FV} in the flat case) and derive that $w$ satisfies, in a system of local coordinates $(\gs_1,...,\gs_{N-1})$ near $\gs_0$,
  $$\BA {lll}\displaystyle
  \CL w:=-\sum_{\ell,j}\myfrac{\prt }{\prt \gs_\ell}\left(A_{j,\ell}\myfrac{\prt w}{\prt \gs_j}\right)+
  \sum_{j}C_j\myfrac{\prt w}{\prt \gs_\ell}+Cw\geq 0,
  \EA$$
  where the matrix $(A_{j,\ell})$ is smooth, symmetric and positive near $\gs_0$ and the $C_j$ and $C$ are bounded. Hence $w$ is locally zero. By a standard argument of  connectedness, this implies that the zero set of $w$ must be empty, contradiction. Hence $\gb=\gb'$.\qeda\medskip
  
   \subsection{Uniqueness of eigenfunction}
   
     The proof is based upon a delicate adaptation of the characterisation of the $p$-Martin boundary obtained  in \cite{lip}, but we first give a proof in the convex case.
   
   \subsubsection{The convex case}

%%%%CONVEX%%%%%%%%%%%%%%%%%%%%%%%%%%%%%%%%%%%%%%%%%%%%%%%%%%%%%%%%%%%%%%%%%%%%%%%%%%%%%%%%%%%%%%%%%%%%%%%%%%%%%%
\bth{convex} Assume $S$ is a convex spherical subdomain. Then two positive solutions of $(\ref{I-1})$ are proportional.
\es
\nind\Proof We recall that a domain $S$ is (geodesically) convex if a minimal geodesic joining two points of $S$ is contained in $S$. If $S\subset S^{N-1}$ is convex, the cone $C_S$ is convex too. Since $S$ is convex, it is Lipschitz and by Theorem D, $\gb_S=\hat\gb_S:=\gb$. Let $\gw$ and $\gw'$ be two positive solutions of $(\ref{I-1})$ satisfying $\sup_S\gw=\sup_S\gw'=1$. We denote by $u_\gw(x)=|x|^{-\gb}\gw(.)$ and $u_{\gw'}(x)=|x|^{-\gb}\gw'(.)$ the corresponding separable $p$-harmonic functions defined in $C_S$. If $0<a<b$, we set $C_S^{a,b}=C_S\cap (B_b\setminus \overline B_a)$. Then for $0<\ge<1$ we denote by $u_\ge$ the unique function which satisfies 
\bel{V-1c}\BA{ll}
-\Gd_pu_\ge=0&\qquad\text{in } C_S^{\ge,1}
\\\phantom{-\Gd_p}
u_\ge=\ge^{-\gb}\gw&\qquad\text{in }  C_S\cap \prt B_\ge
\\\phantom{-\Gd_p}
u_\ge=0&\qquad\text{in }  \left(C_S\cap \prt B_1\right)\cup\left(\prt C_S\cap (\overline B_1\setminus   B_\ge)\right).
  \EA\ee
Then 
\bel{V-2c}(u_\gw-1)_+\leq u_\ge\leq u_\gw\qquad\text{in } C_S^{\ge,1}. 
\ee
Furthermore $\ge\mapsto u_\ge$ is increasing. When $\ge\downarrow 0$, $u_\ge\uparrow u_0$ where $u_0$ is positive and  $p$-harmonic in $C_S^{1,0}$, vanishes on $\prt C_S^{1,0}\setminus\{0\}$ and satisfies $(\ref{V-1c})$ with $\ge=0$. In particular
\bel{V-3c}\displaystyle
\lim_{r\to 0}r^{\gb}u_0(r,\gs)=\gw(\gs)\qquad\text{locally uniformly in  }S.
\ee
We construct the same approximation $u'_\ge$ in $C_S^{\ge,1}$ with $\gw'$  instead of $\gw$.
Mutadis mutandis $(\ref{V-2c})$ holds and $u'_\ge\uparrow u'_0$ which is positive and $p$-harmonic in $C_S^1$, satisfies 
$$(u_{\gw'}-1)_+\leq u'_0\leq u_{\gw'}\qquad\text{in } C_S^{1,0},
$$
and thus
\bel{V-4c}\displaystyle
\lim_{r\to 0}r^{\gb}u'_0(r,\gs)=\gw'(\gs)\qquad\text{locally uniformly in  }S.
\ee
However, by \cite[Th 4]{lip} $u_0$ and $u'_0$ are proportional. Combined with $(\ref{V-3c})$, $(\ref{V-4c})$ it implies the claim.

\subsubsection{Proof of Theorem E}
In what follows we borrow most of our construction from \cite{lip} that we adapt to the case of an infinite cone a make explicit for the sake of completeness. The next {\it nondegeneracy property} of positive $p$-harmonic functions is proved in \cite[Lemma 4.28]{lip}. %%%%%%%%%%%%%%%%%%%%%%%%%%%%%%%%%%%%%%%%%%%%%%%%%%%%%%%
\bprop{nondeg} Let $\Gw\subset\BBR^N$ be a bounded Lipschitz domain and $1<p<\infty$. Then there exist constants $\gr>0$, $c_{14},c_{15}>0$ depending respectively on $\Gw$ (for $\gr$), and $p$, $N$ and the Lipschitz norm $M$ of $\prt \Gw$ (for $c_{14}$ and $c_{15}$) with the property that for any $w\in\prt \Gw$ and any positive $p$-harmonic function $u$ in $\Gw$, continuous in $\overline\Gw\cap \overline B_{2\gr}(w)$ and vanishing on $\prt\Gw\cap B_\gr(w)$, one can find $\xi\in S^{N-1}$, independent of $u$, such that
\bel{V-2}\displaystyle
c_{14}^{-1}\myfrac{u(y)}{\dist\left(y,\prt\Gw\right)}\leq\langle\nabla u(y),\xi\rangle\leq \abs{\nabla u(y)}\leq c_{14}\myfrac{u(y)}{\dist\left(y,\prt\Gw\right)},
\ee
for all $y\in  C_S\cap \overline B_{\frac{\gr\abs w}{c_{15}}}(w)$.
\es

If $\Gw$ is replaced by a cone $C_S$, the nondegeneracy property still holds uniformly on $\prt C_S\setminus\{0\}$. 

\bcor{nondeg-1} Let  $1<p<\infty$, $S\subset S^{N-1}$ is a Lipschitz domain and $C_S$ the cone generated by $S$. \smallskip

\nind (i) Then there exist constants $\gr<\frac12$, $c_{14},c_{15}>0$ depending respectively on $S$  (for $\gr$), and $p$, $N$ and the Lipschitz norm $M$ of $\prt S$ and $diam (S)$ (for $c_{14}$ and $c_{15}$) with the property that for any $w\in\prt C_S$ and any positive $p$-harmonic function $u$ in $C_S$, continuous in $\overline C_S\cap \overline B_{2\gr\abs w}(w)$ and vanishing on $\prt C_S\cap \overline B_{\gr\abs w}(w)$ continuous, one can find $\xi\in S^{N-1}$, independent of $u$, such that
\bel{V-2-a}\displaystyle
c_{14}^{-1}\myfrac{u(y)}{\dist\left(y,\prt C_S\right)}\leq\langle\nabla u(y),\xi\rangle\leq \abs{\nabla u(y)}\leq c_{14}\myfrac{u(y)}{\dist\left(y,\prt C_S\right)},
\ee
for all $y\in B_{\frac{\gr}{c_{15}}}(w)\cap C_S$.\smallskip

\nind (ii) Then there exist positive constants $\gk$ and $c_{16},c_{17}$ depending on $S$  (for $\gk$), and $p$, $N$ and the Lipschitz norm $M$ of $\prt S$ and $diam (S)$ (for $c_{16}, c_{17}$ such that for any $a>0$ and any positive $p$-harmonic function $u$ in $C^a_S$  vanishing on $\prt C_S\cap B_a^c$, there holds
\bel{V-2-b}\displaystyle
c_{16}^{-1}\myfrac{u(y)}{\dist\left(y,\prt C_S\right)}\leq \abs{\nabla u(y)}\leq c_{16}\myfrac{u(y)}{\dist\left(y,\prt C_S\right)}
\qquad\forall y\in C^{c_{17}a}_S\;\text { s.t. }\dist\left(y,\prt C_S\right)\leq \gk\abs y.
\ee
\es

Let $\gw,\gw'\in W^{1,p}_0(S)\cap C(\overline S)$ be positive solutions $(\ref{I-1})$. Since $\frac{\gw}{\gw'}$ is bounded from above and from below in $S$ by positive constants, we can assume, as in the proof of 
Theorem D, that $\gw\geq\gw'$ in $S$ and that the graphs of $\gw$ and $\gw'$ are tangent. hence, if $\gw\neq\gw'$, then $\gw>\gw'$ in $S$ and there exists a sequence $\{\gs_n\}$ converging to $\gs_0\in\prt S$ as $n\to\infty$ such that 
$$\displaystyle \lim_{n\to\infty}\myfrac{\gw'(\gs_n)}{\gw(\gs_n)}=1.
$$
We define $\gd_1=\sup\{\gd>0:\gd\gw<\gw'\}$. For $t\in (\gd_1,1)$, we set
\bel{V-2-1}\displaystyle
\phi_t=\sup \left\{\gw',t\gw\right\}\quad\text{and }\;\psi_t=\inf\left\{\frac{t}{\gd_1}\gw',\gw\right\}
\ee
We also set
\bel{V-2-2}\displaystyle
v_{\phi_t}(r,\gs)=r^{-\gb}\phi_t(\gs)\quad\text{and }\; v_{\psi_t}(r,\gs)=r^{-\gb}\psi_t(\gs)\qquad\forall\,(r,\gs)\in (0,\infty)\ti S.
\ee
\blemma{supinf} The functions $\phi_t$ and $\psi_t$ are respectively a subsolution and a supersolution of $(\ref{I-1})$ in 
$W^{1,p}_0(S)$,  $v_{\phi_t}$ and $v_{\psi_t}$ are respectively a subsolution and a supersolution of $-\Gd_p$ in $C_S$, and there exists $\eta\in W^{1,p}_0(S)$ solution of $(\ref{I-1})$ such that
\bel{V-2-3}\displaystyle
\gw'\leq\phi_t\leq\eta\leq \psi_t\leq\gw\qquad\forall\, t\in (\gd_1,1).
\ee
If $S_t$ is the subset of $\eta\in W^{1,p}_0(S)$ solutions of $(\ref{I-1})$ and satisfying ($\ref{V-2-3})$, then $\gw_t=\sup\{\eta:\eta\in  S_t\}$ belongs to $S_t$. It is increasing with respect to $t$ with uniform limits $\gw'$ when $t\downarrow\gd_1$ and $\gw$ when $t\uparrow 1$. Finally, if $\gth_t=\frac{t-\gd_1}{1-\gd_1}$, there holds
\bel{V-2-4}\displaystyle
\phi_t\leq\gth_t\gw+(1-\gth_t)\gw'\leq \psi_t.
\ee
\es

\nind\Proof Clearly $\phi_t$ and $\psi_t$ are respectively a subsolution and a supersolution of the operator $\CT$, they belong to $W^{1,p}_0(S)\cap L^\infty(S)$ and they satisfy
$\gw'\leq\phi_t\leq \psi_t\leq\gw$. Furthermore, by Dini convergence theorem 
$$ \lim_{t\uparrow 1}\phi_t=\gw=\lim_{t\uparrow 1}\psi_t\quad\text{and }\;\lim_{t\downarrow \gd_1}\phi_t=\gw'=\lim_{t\downarrow \gd_1}\psi_t,
$$
uniformly in $\overline S$. Moreover, in spherical coordinates,
$$\BA {lll}-\Gd_pu(r,\gs)=\left(\left(u_r^2+r^{-2}\abs{\nabla' u}^2\right)^{\frac{p-2}{2}}u_r\right)_r-\myfrac{N-1}{r}\left(u_r^2+r^{-2}\abs{\nabla' u}^2\right)^{\frac{p-2}{2}}u_r\\[4mm]
\phantom{-------------}
-\myfrac{1}{r^2}div'\left(\left(u_r^2+r^{-2}\abs{\nabla' u}^2\right)^{\frac{p-2}{2}}\nabla' u\right).
\EA$$
Hence, if $u(r,\gs)=r^{-\gb}\eta(\gs)$, 
$$-\Gd_pu(r,\gs)=\gb^{p-2}r^{-(p-1)(\gb+1)-1}\CT(\eta).
$$
Thus $v_{\phi_t}$ is a subsolution $-\Gd_p$ in $C_S$ and $v_{\psi_t}$ is a supersolution. 
Since the operator $\CT$ is a Leray-Lions operator, it follows by \cite{BMP} that there exists 
$\eta\in W^{1,p}_0(S)\cap L^\infty (S)$ satisfying $\CT(\eta)=0$ and $\phi_t\leq\eta\leq \psi_t$ in $S$. We denote by 
$S_t$ the set of $\eta\in W^{1,p}_0(S)\cap L^\infty (S)$ satisfying $\CT(\eta)=0$ and $\phi_t\leq\eta\leq \psi_t$ in $S$. Then there exists a sequence $\{\eta_n\} \subset S_t$ and $\gw_t\in W^{1,p}_0(S)\cap L^\infty (S)$ such that $\eta_n(\gs)\uparrow\gw_t(\gs)$ for all $\gs\in \Gs$, where $\Gs$ is a countable dense subset of $S$. By \rlemma {est} $\{\eta_n\}$ is bounded in 
$L^p(S)$, hence in $C^\gg(S)$ for some $\gg\in (0,1)$. By the estimates of the proof of Theorem B-Step 2, $\{\eta_n\}$ is bounded in $W^{1,p}_0(S)$. By standard regularity theory, we can also assume that $\eta_n\to\gw_t$ in the $C^1_{loc}(S)$-topology. Hence $\gw_t$ is a weak solution of $(\ref{I-1})$, it belongs to $W^{1,p}_0(S)\cap L^\infty (S)$  and satisfies 
$\phi_t\leq\gw_t\leq \psi_t$. Therefore it is the maximal element of $S_t$. The monotonity of $\gw_t$ is a consequence of the monotonicity of $\phi_t$ and $\psi_t$ and the last statement $(\ref{V-2-4})$  is a straightforward computation. \qeda\medskip

Next we recall the deformation of $p$-harmonic functions already used in \cite{lip}. If $\gt\in (0,1)$ and $0<a<b$, we denote by $v_{\gt,a,b}$ the $p$-harmonic function defined in $C_S^{a,b}$ satisfying 
\bel{V-2-5}\BA {lll}\displaystyle
v_{\gt,a,b}(x)=\left\{\BA {lll}a^{-\gb}(\gt\gw+(1-\gt)\gw')(\frac{x}{\abs x})\quad &\text{if }x\in C_S\cap \prt B_a\\[2mm]
0\quad &\text{if }x\in C_S\cap \prt B_b\\[2mm]
0&\text{if }x\in \prt C_S\cap \left(\overline B_b\setminus B_a\right). 
\EA\right.\EA\ee
%%%%%%%%%%%%%%%%%%%%%%%%%%%%%%%%%%%%%%%%%%%%%%%%%%%%%%%%%%%%%%%%%%%%%%%%%%%%%%%%%%%%%%%%%%%%%%%%%%%%%%%%%%%%%

\blemma{deform} The mapping $(\gt,b)\mapsto v_{\gt,a,b}$ is continuous and increasing. If $v_{\gt,a}=\displaystyle\lim_{b\to\infty} v_{\gt,a,b}$, then it is a positive p-harmonic function in $C_S^{a,\infty}$ vanishing on $\prt S\cap B_a^c$, and
there holds
\bel{V-2-7}\BA {lll}
u_{\gw'}(x)\leq v_{\phi_{\gt^*}}(x)\leq v_{\gt,a}(x)\leq v_{\psi_{\gt^*}}(x)\leq u_\gw(x)\quad\forall\, x\in C_S^{a,\infty},
\EA\ee 
where $\gt^*=(1-\gd_1)\gt+\gd_1$ and as a consequence  
\bel{V-2-6}\BA {lll}\displaystyle
\lim_{\gt\uparrow 1}\sup_{\abs x\geq a}\abs x^{\gb}\left(u_\gw(x)-v_{\gt,a}(x)\right)=0\quad\text{and }\;
\lim_{\gt\downarrow 0}\sup_{\abs x\geq a}\abs x^{\gb}\left(v_{\gt,a}(x)-u_{\gw'}(x)\right)=0
\EA\ee 
Furthermore  
\bel{V-2-8}\BA {lll}\displaystyle
0\leq \myfrac{v_{\gt',a}-v_{\gt,a}}{\gt'-\gt}\leq \left(\myfrac{1}{\gd_1}-1\right)v_{\gt',a}\quad\forall\, 0\leq \gt<\gt'\leq 1.
\EA\ee
\es

\nind\Proof The uniqueness and the (strict) monotonicity  of $(\gt,b)\mapsto v_{\gt,a,b}$ follow from the monotonicity of $\gt\mapsto\gt\gw+(1-\gt)\gw' $ and the strong maximum principle. The continuity is a consequence of uniqueness and regularity theory for $p$-harmonic functions. It follows from $(\ref{V-2-4})$ with $t=\gt^*$ and the fact that $v_{\phi_{\gt^*}}$ and 
$v_{\psi_{\gt^*}}$ are respectively a subsolution and a supersolution of $-\Gd_p$, that we have
$$u_{\gw'}(x)\leq v_{\phi_{\gt^*}}(x)\leq v_{\gt,a,b}(x)\leq v_{\psi_{\gt^*}}(x)\leq u_\gw(x)\quad\forall\, x\in C_S^{a,b},$$
which yields $(\ref{V-2-7})$. Similarly, we have on $\prt C_S^{a,b}$ 
\bel{V-2-9}0\leq\myfrac{v_{\gt',a,b}-v_{\gt,a,b}}{\gt'-\gt}=u_{\gw}-u_{\gw'}\leq (\gd_1^{-1}-1)u_{\gw'}\leq (\gd_1^{-1}-1)v_{\gt,a,b},
\ee
equivalently
\bel{V-2-10}0\leq v_{\gt',a,b}\leq
\left(1+(\gt'-\gt) (\gd_1^{-1}-1)\right)v_{\gt,a,b}.
\ee
By the maximum principle $(\ref{V-2-9})$ holds in $C_S^{a,b}$. This implies $(\ref{V-2-8})$.\qeda\\

As a consequence of $(\ref{V-2-8})$, $\prt_\gt v_{\gt,a}$ exists for almost all $\gt\in (0,1)$ in $W^{1,p}_0(C_S^{a,b})$ for all $b>a$ and it is a solution of  
\bel{V-3}\BA {lll}\displaystyle
\BBL w=\nabla.\left((p-2)\abs{\nabla v_{\gt,a}}^{p-4}\langle\nabla v_{\gt,a}.\nabla Z\rangle\nabla v_{\gt,a}\right)
\\[4mm]\phantom{\CL w}\displaystyle
=\sum_{i,j}\myfrac{\prt}{\prt x_j}\left(b_{i,j}(x)\myfrac{\prt w}{\prt x_i}\right)=0
\EA\ee
where
$$b_{i,j}(x)=\abs{\nabla v_{\gt,a}}^{p-4}\left((p-2)\myfrac{\prt v_{\gt,a}}{\prt x_j}\myfrac{\prt v_{\gt,a}}{\prt x_i}+\gd_{ij}\abs{\nabla v_{\gt,a}}^{2}\right).
$$
$\BBL $ satisfies the following ellipticity condition
\bel{V-4}\BA {lll}\displaystyle
 \min\{1,p-1\}\abs{\nabla v_{\gt,a}}^{2}\abs\xi^2\leq\sum_{i,j}b_{i,j}(x)\xi_i\xi_j\leq \max\{1,p-1\}\abs{\nabla v_{\gt,a}}^{2}\abs\xi^2\quad\forall \xi\in\BBR^N.
\EA\ee
It is important to notice that $\BBL v_{\gt,a}=(p-1)\Gd_pv_{\gt,a}=0$. The estimate $(\ref{V-4})$ combined with $(\ref{V-2-b})$ 
and the decay of $v_{\gt,a}$ and $\prt_\gt v_{\gt,a}$ implies that they satisfy Harnack inequality and boundary Harnack inequality in $C_S^a$. There exists a constant $\hat c>c_{17}>1$ (see $\ref{V-2-b}$) such that 
\bel{V-5}\BA {lll}\displaystyle
\myfrac{1}{\hat c}\myfrac{\prt_\gt v_{\gt,a}(x_a)}{v_{\gt,a}(x_a)}\leq 
\myfrac{\prt_\gt v_{\gt,a}(x)}{v_{\gt,a}(x)}\leq \hat c\myfrac{\prt_\gt v_{\gt,a}(x_a)}{v_{\gt,a}(x_a)}\qquad\forall x\in C_S^{\hat ca},
\EA\ee
where $x_a=(\hat ca,\gs_0)$ for some $\gs_0\in S$ fixed. We set 
\bel{V-6}\BA {lll}\displaystyle
M(t)=\sup_{x\in C^t_S}\myfrac{\prt_\gt v_{\gt,a}(x)}{v_{\gt,a}(x)}\quad\text{and }\;\;m(t)=\inf_{x\in C^t_S}\myfrac{\prt_\gt v_{\gt,a}(x)}{v_{\gt,a}(x)}\qquad\forall t>a
\EA\ee
\blemma{decay} For $t>\hat ca$ there holds
\bel{V-7}\BA {lll}\displaystyle
M(\hat c t)-m(\hat c t)\leq \myfrac{\hat c^2-1}{\hat c^2+1}\left(M(t)-m(t)\right).
\EA\ee
\es
\nind\Proof There holds 
$$\prt_\gt v_{\gt,a}-m(t)v_{\gt,a}\geq 0\quad\text{and }\;\;M(t)v_{\gt,a}-\prt_\gt v_{\gt,a}\geq 0\qquad\in C^t_S.
$$
Estimate $(\ref{V-5})$ is valid for any couple of positive solutions $(h_1,h_2)$ of $\BBL h=0$ in $C_S^a$ vanishing on $\prt C_S^a\cap B_a^c$, in particular for $\left(\prt_\gt v_{\gt,a}-m(t)v_{\gt,a},v_{\gt,a}\right)$ and 
$\left(M(t)v_{\gt,a}-\prt_\gt v_{\gt,a},v_{\gt,a}\right)$. Hence 
\bel{V-7}\BA {lll}
\myfrac{1}{\hat c}\left(\myfrac{\prt_\gt v_{\gt,a}(x_a)}{v_{\gt,a}(x_a)}-m(t)\right)\leq \myfrac{\prt_\gt v_{\gt,a}(x)}{v_{\gt,a}(x)}-m(t)\leq 
\hat c\left(\myfrac{\prt_\gt v_{\gt,a}(x_a)}{v_{\gt,a}(x_a)}-m(t)\right)\quad\forall x\in C^t_S.
\EA\ee
This implies
$$\myfrac{1}{\hat c}\left(\myfrac{\prt_\gt v_{\gt,a}(x_a)}{v_{\gt,a}(x_a)}-m(t)\right)\leq m(\hat ct)-m(t),
$$
and
$$\myfrac{\prt_\gt v_{\gt,a}(x)}{v_{\gt,a}(x)}-m(t)\leq \hat c^2(m(\hat ct)-m(t))\quad\forall x\in C^t_S.
$$
Finally
\bel{V-8}\BA {lll}
M(\hat ct)-m(t)\leq \hat c^2(m(\hat ct)-m(t)). 
\EA\ee
Similarly
\bel{V-9}\BA {lll}
M(t)-m(\hat ct)\leq \hat c^2(M(t)-M(\hat ct)).
\EA\ee
Summing the two inequalities we get
$$\left(M(t)-m(t)\right)+\left(M(\hat ct)-m(\hat ct)\right)\leq \hat c^2\left(\left(M(t)-m(t)\right)-\left(M(\hat ct)-m(\hat ct)\right)\right),
$$
which yields $(\ref{V-7})$.\qeda
\medskip

\nind{\it End of the proof.} By the differentiability property of $v_{\gt,a}$ with respect to $\gt$,  there exists two countable dense sets $\{(r_\gn\}\subset [a,\infty)$ and $\{\gs_\gm\}\subset [a,\infty)$ such that $\prt_\gt v_{\gt,a}(r_\gn,\gs_\gm)$ exists for almost all $\gt$. We put $x_{\gn,\gm}=(r_\gn,\gs_\gm)$, hence
\bel{V-10}\BA {lll}
\ln\left(\myfrac{\gw(\gs_\gm)}{\gw'(\gs_\gm)}\right)-\ln\left(\myfrac{\gw(\gs_{\gm'})}{\gw'(\gs_{\gm'})}\right)
=
\ln\left(\myfrac{v_{1,a}(x_{\gn,\gm})}{v_{0,a}(x_{\gn,\gm})}\right)-\ln\left(\myfrac{v_{1,a}(x_{\gn,\gm'})}{v_{0,a}(x_{\gn,\gm'})}\right)\\[4mm]
\phantom{\ln\left(\myfrac{\gw(\gs_\gm)}{\gw'(\gs_\gm)}\right)-\ln\left(\myfrac{\gw(\gs_{\gm'})}{\gw'(\gs_{\gm'})}\right)}
=\myint{0}{1}\left(
\myfrac{\prt_\gt v_{\gt,a}(x_{\gn,\gm})}{v_{\gt,a}(x_{\gn,\gm})}-\myfrac{\prt_\gt v_{\gt,a}(x_{\gn,\gm'})}{v_{\gt,a}(x_{\gn,\gm'})}\right)d\gt.
\EA\ee
Using the continuity of $\frac{\gw}{\gw'}$ and the density of $\{\gs_m\}$ we derive
\bel{V-11}\BA {lll}\abs{
\ln\left(\myfrac{\gw(\gs)}{\gw'(\gs)}\right)-\ln\left(\myfrac{\gw(\gs')}{\gw'(\gs')}\right)}.
\leq M(r_\gn)-m(r_\gn)\qquad \forall(\gs,\gs')\in S\ti S.
\EA\ee
We can assume that $r_\gn\geq \hat c^{\gn_n}a$ for some sequence $\{\gn_n\}$ tending to infinity with $n$, hence
\bel{V-12}\BA {lll}\abs{
\ln\left(\myfrac{\gw(\gs)}{\gw'(\gs)}\right)-\ln\left(\myfrac{\gw(\gs')}{\gw'(\gs')}\right)}
\leq \gth^n\left(M(\hat c^{\gn_1})-m(\hat c^{\gn_1})\right)\qquad \forall(\gs,\gs')\in S\ti S\quad\forall n\in\BBN^*,
\EA\ee
where $\gth=\frac{\hat c^2-1}{\hat c^2+1}<1$. Letting $n\to\infty$ implies the claim.\qeda
%%%%%%%%%%%%%%%%%%%%%%%%%%%%%%%%%%%%%%%%%%%%%%%%%%%%%%%%%%%%%%%%%%%%%%%%%%%%%%%%%%%%%%%%%%%%%%%%%%%%%%%%%%%%%%%%%%%%%%%%%%%%%%%%%%%%%%%%%%%%%%%%%%%%%%%%%%%%%%%%%%%%%%%%%%%%%%%%%%%JUNK%%%%%%%%%%%%%%%%%%%%%%%%%%%%%%%%%%%%%%%%%%%%%%%%%%%%%%%%%%%%%%%%%%%%%%%%%%%%%%%%%%%%%%%%%%%%%%%%%%%%%%%%%%%%%%%%%%%%%%%%%%%%%%%%%%%%%%%%%%%%%%%%%%%%%
 
 %%%%%%%%%%%%%%%%%%%%%%%%%%%%%%%%%%%%%%%%%%%%%%%%%%%%%%%%%%%%%%%%%%%%%%%%%%%%%%%%%%%%%%%%%%%%%%%%%%%%%%%%%%%%%%%%%%%%%%%%%%%%%%%%%%%%%%%%%%%%%%%%%%%%%%%%%%%%%%%%%%%%%%%%%%%%%%%%%%%%%%%%%%%%%%%%%%%%%%%%%%%%%%%%%%%%%%%%%%%%

 %%%%%%%%%%%%%%%%%%%%%%%%%%%%%%%%%%%%%%%%%%%%%%%%%%%%%
%%%%%%%%%%%%%%%%%%%%%%%%%%%%%%%%%%%%%%%%%%%%%%%%%%%%%
 %%%%%%%%%%%%%%%%%%%%%%%%%%%%%%%%%%%%%%%%%%%%%%%%%%%%%
%%%%%%%%%%%%%%%%%%%%%%%%%%%%%%%%%%%%%%%%%%%%%%%%%%%%%
 %%%%%%%%%%%%%%%%%%%%%%%%%%%%%%%%%%%%%%%%%%%%%%%%%%%%%
%%%%%%%%%%%%%%%%%%%%%%%%%%%%%%%%%%%%%%%%%%%%%%%%%%%%%

  %%%%%%%%%%%%%%%%%%%%%%%%%%%%%%%%%%%%%%%%%%%%%%%%%%%%%
%%%%%%%%%%%%%%%%%%%%%%%%%%%%%%%%%%%%%%%%%%%%%%%%%%%%%
%%%%%%%%%%%%%%\blemma{est}%%%%%%%%%%%%%%%%%%%%%%%%%%%%%%%%%%%%%%%
%%%%%%%%%%%%%%%%%%%%%%%%%%%%%%%%%%%%%%%%%%%%%%%%%%%%%
%%%%%%%%%%%%%%%%%%%%%%%%%%%%%%%%%%%%%%%%%%%%%%%%%%%%%

\nind G. K.: Centro de Modelamiento Matem\`atico (UMI 2807 CNRS), Universidad de Chile, Casilla 170 Correo 3, Santiago, Chile.
\medskip

\nind L. V. : Laboratoire de Math\'ematiques et Physique Th\'eorique, (UMR 7350 CNRS), Facult\'e des Sciences, Universit\'e Fran\c{c}ois Rabelais, 37200, Tours, France. 

\end{document}